\documentstyle[draft%
,amscd%
,syntonly%
,amssymb%
]{amsart}

\newtheorem{thm}{Theorem}   

\newtheorem{lem}[thm]{Lemma}

\newtheorem{rem}[thm]{Remark}
\newtheorem{defn}{Definition}

\begin{document}


\title[Tunnel Number] {The Tunnel Number of the Sum of
$\bf \Huge n$  knots is at least $\bf \Huge n$} 
\author{Martin Scharlemann \\
Jennifer Schultens} 
\address{Department of Mathematics \\ 
UCSB \\
Santa Barbara, CA 93106\\ 
$\hspace{32 mm}$ Department of Mathematics
\& CS \\ 
Emory University \\ 
Atlanta, GA 30322}
\email{mgscharl@@math.ucsb.edu \\ 
jcs@@mathcs.emory.edu} 

\thanks{The first author is supported by an NSF grant. \\ 
The second author is supported by an NSF postdoctoral fellowship.}  
\date{\today}

\begin{abstract}
We prove that the tunnel number of the sum of $n$ knots is
at least $n$.
\end{abstract}

\maketitle

\section{Introduction}

\vspace{5 mm}

In \cite{N}, Norwood showed that tunnel number $1$ knots are prime.
This led to the more general conjecture, see for instance
\cite[Problem 1.70B]{K}, that the tunnel number of a sum of $n$
knots is at least $n$.  Here we prove this conjecture.  The idea
is to show that the splitting surface of a Heegaard splitting
corresponding to a tunnel system realizing the tunnel number of the 
sum of $n$ knots intersects each individual knot complement
essentially.  Then a sophisticated Euler characteristic argument,
based on the idea of untelescoping the Heegaard splitting, yields the
result.  

We wish to thank MSRI, where part of this work was carried out.

\vspace{5 mm}

\section{Preliminaries}

\vspace{5 mm}

For standard definitions concerning knots, see \cite{BZ} or \cite{R}
and for those concerning $3$-manifolds, see \cite{H} or \cite{J}.

\begin{defn}
Let $N$ be a submanifold of $M$, we denote an open regular
neighborhood of $N$ in $M$ by $\eta(N)$.
\end{defn}

\begin{defn}
Let $K$ be a knot in $S^3$.  Denote the \underline{complement of K},
$S^3$ $-$ $\eta(K)$, by $C(K)$.  
\end{defn}

\begin{rem}  \label{rem:dec}
Let $K$ $=$ $K_1$ $\#$ $K_2$ be the sum of two knots.  Then
the decomposing sphere gives rise to a decomposing annulus $A$
properly embedded in $C(K)$ such that $C(K)$ $=$ $C(K_1)$ $\cup_A$
$C(K_2)$.  If $K$ $=$ $K_1$ $\#$ $\dots$ $\#$ $K_n$, then we may assume
that the decomposing spheres are nested, so that $C(K)$ $=$ $C(K_1)$
$\cup_{A_1}$ $\dots$ $\cup_{A_{n-1}}$ $C(K_n)$.  
\end{rem}

\begin{defn} \label{defn:tunnel system}
A \underline{tunnel system} for a knot $K$ is a collection of disjoint
arcs $\cal T$ $=$ $t_1$ $\cup$ $\dots$ $\cup$ $t_n$ properly embedded
in $C(K)$ such that $C(K)$ $-$ $\eta(\cal T)$ is a handlebody.  The
\underline{tunnel number of $K$}, denoted by $t(K)$, is the least
number of arcs required in a tunnel system for $K$.   
\end{defn}
 
\begin{defn} \label{defn:cb}
A \underline{compression body} is a $3$-manifold $W$ obtained from a
connected closed orientable surface $S$ by attaching $2$-handles to
$S$ $\times$ $\{0\}$ $\subset$ $S$ $\times$ $I$ and capping off any
resulting $2$-sphere boundary components.  We denote
$S$ $\times$ $\{1\}$ by $\partial_+W$ and $\partial W$ $-$ $\partial_+W$
by $\partial_-W$. 
\end{defn}

\begin{defn}
A \underline{set of defining disks} for a compression body $W$ is a
set of disks $\{D_1, \dots, D_n\}$ properly embedded in $W$ with
$\partial D_i$ $\subset$ $\partial_+W$ for $i$ $=$ $1$, $\dots,$ $n$
such that the result of cutting $W$ along $D_1$ $\cup$ $\dots$ $\cup$
$D_n$ is homeomorphic to $\partial_-W$ $\times$ $I$.
\end{defn}

\begin{defn} \label{defn:Heegaard splitting}
A \underline{Heegaard splitting} of a $3$-manifold $M$ is a
decomposition $M$ $=$ $V$ $\cup_S$ $W$ in which $V$, $W$ are
compression bodies such that $V$ $\cap$ $W$ $=$ $\partial_+V$ $=$
$\partial_+W$ $=$ $S$ and $M$ $=$ $V$ $\cup$ $W$.  We call $S$ the
\underline{splitting surface} or \underline{Heegaard surface}.
\end{defn}

\begin{defn} 
Let $M$ $=$ $V$ $\cup_S$ $W$ be an irreducible Heegaard splitting.  We
may think of $M$ as being obtained from $\partial_-V$ $\times$ $I$ by
attaching all $1$-handles dual to $2$-handles in $V$ followed by all
$2$-handles in $W$, followed, perhaps, by $3$-handles.  An
\underline{untelescoping} of $M$ $=$ $V$ $\cup_S$ $W$ is a
rearrangement of the order in which the $1$-handles (of $V$) and the
$2$-handles (dual to the $1$-handles of $W$) are attached.  This
rearrangement is chosen so that $M$ is decomposed into submanifolds
$M_1$, $\dots$, $M_m$, such that $M_i$ $\cap$ $M_{i+1}$ $=$ $F_i$ and
$F_i$ is an incompressible surface in $M$, and such that the $M_i$
inherit, from a subcollection of the original $1$-handles and
$2$-handles, strongly irreducible Heegaard splittings $M_1$ $=$ $V_1$
$\cup_{S_1}$ $W_1$, $\dots,$ $M_m$ $=$ $V_m$ $\cup_{S_m}$ $W_m$.
Unless $M$ is a lens space or $S^1$ $\times$ $S^2$, no $S_1,$ $\dots,$
$S_m$ is a torus.  For details see \cite{ST1} and \cite{S1}.  We
denote the untelescoping of $M$ $=$ $V$ $\cup_S$ $W$ by $M$ $=$ $(V_1$
$\cup_{S_1}$ $W_1)$ $\cup_{F_1}$ $\dots$ $\cup_{F_{m-1}}$ $(V_m$
$\cup_{S_m}$ $W_m)$.  For convenience, we will occasionally denote
$\partial_-V$ $=$ $\partial_-V_1$ by $F_0$.
\end{defn}

\begin{lem}  \label{lem:cut}
$\chi(S)$ $=$ $\sum_{i=1}^m$ $\chi(S_i)$ $-$ $\sum_{i=1}^{m-1}$
$\chi(F_i)$.
\end{lem}

\begin{pf}
Let $M$ $=$ $V$ $\cup_S$ $W$ be a Heegaard splitting, then \[\chi(S) =
\chi(\partial_-V) - 2(\#(1-handles\; attached \;in\; V)-\#(0-handles\;
attached \;in\; V))\] and in an untelescoping, \[\chi(S_i) =\] \[
\chi(\partial_-V_i) - 2(\#(1-handles\; attached\; in\;
V_i)-\#(0-handles\; attached \;in\; V_i)) = \]
\[\chi(F_{i-1}) - 2(\#(1-handles \;attached\; in\; V_i)-\#(0-handles\;
attached \;in\; V_i)).\] So, since $1-handles$ are merely reordered in
an untelescoping,
\[\chi(S) = \] \[\chi(\partial_-V) - 2\sum_{i=1}^m\#((1-handles\; 
attached \;in\; V_i)-\#(0-handles\; attached \;in\; V_i)) = \]
\[\chi(\partial_-V) - \sum_{i=1}^m \chi(F_{i-1}) + \sum_{i=1}^m
\chi(S_i).\] 
\end{pf}

\begin{lem}  \label{lem:essential}
Let $P$ be a properly embedded incompressible surface in an
irreducible $3$-manifold $M$ and let $M$ $=$ $(V_1$ $\cup_{S_1}$
$W_1)$ $\cup_{F_1}$ $\dots$ $\cup_{F_{m-1}}$ $(V_m$ $\cup_{S_m}$
$W_m)$ be an untelescoping of a Heegaard splitting $M$ $=$ $V$
$\cup_S$ $W$.  Then $(\cup_{i=1}^{m-1} F_i)$ $\cup$ $(\cup_{i=1}^m
S_i)$ can be isotoped to intersect $P$ only in curves essential in $P$.
\end{lem}

\begin{rem}
This lemma demonstrates the advantage of working with untelescopings
of Heegaard splittings rather than Heegaard splittings.  It is a deep fact
that the splitting surface of a strongly irreducible Heegaard splitting
can be isotoped to intersect a properly embedded incompressible surface
only in curves essential in this surface.  This fact is proven for 
instance in \cite[Lemma 6]{Sc}.
\end{rem}

\begin{pf}
Here $(\cup_{i=1}^{m-1} F_i)$ may be isotoped to intersect $P$ only in
curves essential in $P$ by a standard innermost disk argument, since
both are incompressible.  Then $P_i$ $=$ $P$ $\cap$ $M_i$ is a
properly embedded incompressible surface in $M_i$.  It follows that
each $S_i$ may be isotoped in $M_i$ to intersect $P_i$ only in curves
essential in $P_i$, by \cite[Lemma 6]{Sc}.  Note that the latter
isotopies fix $(\cup_{i=1}^{m-1} F_i)$.
\end{pf}

\begin{lem}  \label{lem:par}
Let $K$ be a prime knot and let $A$ be an annulus properly embedded
in $C(K)$ such that the components of $\partial A$ are meridians.
Then $A$ is boundary parallel.
\end{lem}

\begin{pf}
In $S^3$, $A$ can be extended to a sphere by adding two meridian
disks.  This sphere intersects $K$ in two points.  Since $K$ is
prime, one side of the sphere contains a single unknotted arc.
\end{pf}

\begin{lem}  \label{lem:cbs}
Let $P$ be an incompressible surface in a compression body $W$.  Then
the result of cutting $W$ along $P$ is a collection of compression
bodies.
\end{lem}

\begin{pf}
This is \cite[Lemma 2]{Sc}.
\end{pf}

\begin{rem}
In the above lemma, $P$ need not be connected.
\end{rem}

\begin{lem}  \label{lem:max}
If $\cal A$ is a collection of incompressible annuli in a
compression body $W$, then in any component $X$ of $W$ $-$ $\cal A$,
$\chi(\partial_+W$ $\cap$ $X)$ $\leq$ $\chi(\partial_-W$ $\cap$ $X)$.
\end{lem}

\begin{pf}
Let $\cal D$ be a set of defining disks for $W$.  We argue by
induction on the pair $(\vline \chi(\partial_-W) - \chi(\partial_+W)
\vline,$ $\vline \cal A \cap \cal D \vline)$.  If $\vline
\chi(\partial_-W) - \chi(\partial_+W) \vline$ $=$ $0$, then ($\cal D$ $=$
$\emptyset$ and) all annuli are spanning annuli and the result follows.

To complete the inductive step, suppose there is a disk $D$ in $\cal
D$ such that $D$ $\cap$ $\cal A$ $=$ $\emptyset$.  The result of
cutting $W$ along $D$ is a compression body $W'$ with $\vline
\chi(\partial_-W') - \chi(\partial_+W') \vline$ $<$ $\vline
\chi(\partial_-W) - \chi(\partial_+W) \vline$, or two compression
bodies $W'$ and $W''$ with $\vline \chi(\partial_-W') -
\chi(\partial_+W') \vline$ $<$ $ \vline \chi(\partial_-W) -
\chi(\partial_+W) \vline$ and $\vline \chi(\partial_-W) -
\chi(\partial_+W'') \vline$ $<$ $\vline \chi(\partial_+W) \vline$.
The components of $W$ $-$ $\cal A$ can be obtained from the components
of $W'$ $-$ $\cal A$ or of $W'$ $-$ $\cal A$ and $W''$ $-$ $\cal A$ by
attaching a $1-handle$ either to a single component or so as to
connect two components.  In both cases, the result follows from the
inductive hypotheses.

If there is no such disk, consider $\cal D$ $\cap$ $\cal A$.  If there
is an arc $\alpha$ in $\cal D$ $\cap$ $\cal A$ that is inessential in
$\cal A$, then we may assume that $\alpha$ is outermost in $\cal A$,
and we may cut the disk $D$ in $\cal D$ containing $\alpha$ along
$\alpha$ and paste on two copies of the disk cut off by $\alpha$ in
$\cal A$ to obtain a new disk $D'$.  Replacing $D$ by $D'$ in $\cal D$
produces a new set of defining disks $\cal D'$ with $\vline \cal A
\cap \cal D' \vline$ $<$ $\vline \cal A \cap \cal D \vline$.

If all arcs in $\cal D$ $\cap$ $\cal A$ are essential in $\cal A$, let
$\beta$ be an arc in $\cal D$ $\cap$ $\cal A$ that is outermost in
$\cal D$.  Let $A$ be the annulus in $\cal A$ that gives rise to $\beta$.
Cutting and pasting $A$ along $\beta$ and the outermost disk cut off
in $\cal D$ yields a disk $D'$ disjoint from $\cal A$.  If $D'$ is
inessential, then $A$ is inessential and can be ignored.  (Since cutting
along $A$ does not alter any components or their Euler characteristics.)
If $D'$ is essential, the result follows as above.  This completes the
inductive step.
\end{pf}

\vspace{5 mm}

\section{The Combinatorics}

\vspace{5 mm}

In the following, we consider a tunnel system $\cal T$, realizing the
tunnel number of $K_1$ $\#$ $\dots$ $\#$ $K_n$.  We also consider the
Heegaard splitting $C(K_1$ $\#$ $\dots$ $\#$ $K_n)$ $=$ $V$ $\cup_S$
$W$ corresponding to $\cal T$ and an untelescoping $C(K_1$ $\#$
$\dots$ $\#$ $K_n)$ $=$ $(V_1$ $\cup_{S_1}$ $W_1)$ $\cup_{F_1}$
$\dots$ $\cup_{F_{m-1}}$ $(V_m$ $\cup_{S_m}$ $W_m)$ of $C(K_1$ $\#$
$\dots$ $\#$ $K_n)$ $=$ $V$ $\cup_S$ $W$.  Set $M_i$ $=$ $V_i$ $\cup$
$W_i$.  By Remark \ref{rem:dec}, $C(K_1$ $\#$ $\dots$ $\#$ $K_n)$ $=$
$C(K_1)$ $\cup_{A_1}$ $\cup$ $\dots$ $\cup_{A_{n-1}}$ $C(K_n)$.  We
will always assume that $\partial_- V_1$ $=$ $\partial C(K_1$ $\#$
$\dots$ $\#$ $K_n)$ and that $\cup_{i=1}^{m-1} F_i$ and $\cup_{i=1}^m
S_i$ intersect $\cup_{j=1}^{n-1} A_j$ only in curves essential in
$\cup_{j=1}^{n-1} A_j$.  We will, furthermore, assume that, subject to
these constraints, the number of intersections of $\cup_{i=1}^{m-1}
F_i$ and $\cup_{i=1}^m S_i$ with $\cup_{j=1}^{n-1} A_j$ is minimal.

\begin{defn}
Set $S_{ij}$ $=$ $S_i$ $\cap$ $C(K_j)$, $F_{ij}$ $=$ $F_i$ $\cap$
$C(K_j)$ and $A_{ij}$ $=$ $M_i$ $\cap$ $A_j$.  
\end{defn}

\begin{lem}  \label{lem:even}
For all $i$,$j$, $\chi(S_{ij})$ and $\chi(F_{ij})$ are even.
\end{lem}

\begin{pf}
Here $F_i$ is separating, so $F_i$ $\cap$ $A_{j-1}$ is separating.
Since $\partial A_{j-1}$ $\subset$ $\partial C(K_1$ $\#$ $\dots$ $\#$
$K_n)$ which is a torus, hence connected, both components lie on one
side of $F_i$, hence $\vline F_i$ $\cap$ $A_{j-1} \vline$ is even.
The same is true for $\vline F_i$ $\cap$ $A_j \vline$.  Thus
$\chi(F_i$ $\cap$ $C(K_j))$ $=$ $2$ $-$ $2(genus(F_i$ $\cap$ $C(K_j)))$
$-$ $\vline F_i$ $\cap$ $(A_{j-1}$ $\cup$ $A_j) \vline$ is even.
Similarly for $S_i$.
\end{pf}

\begin{defn}
Set $x_{ij}$ $=$ $-1/2\chi(F_{ij})$ and
$y_{ij}$ $=$ $-1/2\chi(S_{ij})$.
\end{defn}

\vspace{5 mm}

\begin{center}
\begin{tabular}{|l|l|l|l|l|l|} \hline
  & 1 & $\dots$ & $j$ & $\dots$ & n \\ \hline
1 & & & & & \\ \hline
$\dots$ & & $\dots$ & & $\dots$ & \\ \hline
$i$ & & & $x_{ij}$ & & \\ \hline
$\dots$ & & $\dots$ & & $\dots$ & \\ \hline
$m$-1 & & & & & \\ \hline
\end{tabular}

\vspace{2 mm}
fig. 1a

\vspace{10 mm}

\begin{tabular}{|l|l|l|l|l|l|} \hline
  & 1 & $\dots$ & $j$ & $\dots$ & n \\ \hline
1 & & & & & \\ \hline
$\dots$ & & $\dots$ & & $\dots$ & \\ \hline
$i$ & & & $y_{ij}$ & & \\ \hline
$\dots$ & & $\dots$ & & $\dots$ & \\ \hline
$m$ & & & & & \\ \hline
\end{tabular}

\vspace{2 mm}
fig. 1b

\vspace{5 mm}

\end{center}

\begin{lem}  \label{lem:trap}
Under the assumptions above, $y_{ij}$ $\geq$ $max\{x_{i-1j},
x_{ij}\}$.
\end{lem}

\begin{pf}
This follows from Lemma \ref{lem:max}.
\end{pf}

\begin{lem}  \label{lem:nonzero}
For all $j$, there is an $i$, such that $y_{ij}$ $>$ $0$.
\end{lem}

\begin{pf}
Suppose $y_{ij}$ $=$ $0$ for $i$ $=$ $1,$ $\dots,$ $m$.  Then $x_{ij}$
$=$ $0$ for $i$ $=$ $1,$ $\dots,$ $m-1$.  So \[G_j = (\cup_{i=1}^{m-1}
F_{ij}) \cup (\cup_{i=1}^m S_{ij}) \subset C(K_j)\] is a collection of
annuli and tori.  Since the tori arise only in $\cup_{i=1}^{m-1}
F_{ij}$, they are incompressible separating tori.  Thus if a torus
component $T$ of $F_i$ is in $C(K_j)$, then so is a component of
$S_{i'}$, which cannot be a torus, for some $i'$.  But this would
contradict $y_{i'j}$ $=$ $0$.  Hence $G_j$ consists entirely of
annuli.  By Lemma \ref{lem:par}, the annuli are all
boundary parallel.  Hence cutting $C(K_j)$ along the annular
components of $G_j$ yields a copy of $C(K_j)$.
By Lemma \ref{lem:cbs}, all components of $C(K_j)$ cut along $G$
are compression bodies, a contradiction.  
\end{pf}

\begin{lem}  \label{lem:strict}
For all $j$, \[\sum_{i=1}^{m} y_{ij} > \sum_{i=1}^{m-1} x_{ij}.\]
\end{lem}

\begin{pf}
This follows by comparing the tables in fig. 1.  By Lemma
\ref{lem:trap}, the largest value encountered in a given column of the
table in fig. 1a occurs one time more often in the corresponding column of
the table in fig. 1b.  If the largest value encountered in a column in the
table in fig. 1a is zero, then by Lemma \ref{lem:nonzero}, there must
be nonzero entries in the corresponding column of the table in fig. 1b.
\end{pf}

\begin{rem}
Since all numbers involved are integers, it follows that
$\sum_{i=1}^{m}$ $y_{ij}$ $\geq$ $1$ $+$ $\sum_{i=1}^{m-1}$ $x_{ij}$,
for all $j$.
\end{rem}

\begin{thm}
$t(K_1$ $\#$ $\dots$ $\#$ $K_n)$ $\geq$ $n$.
\end{thm}

\begin{pf}
Here \[\sum_{j=1}^n(\sum_{i=1}^m y_{ij}) \geq
\sum_{j=1}^n(1 + \sum_{i=1}^{m-1} x_{ij}) = n + \sum_{j=1}^n \sum_{i=1}^{m-1}
x_{ij}.\] Hence, \[\sum_{j=1}^n \sum_{i=1}^m y_{ij} -
\sum_{j=1}^n \sum_{i=1}^{m-1} x_{ij} \geq n.\] Thus,
\[\sum_{j=1}^n \sum_{i=1}^m -2(y_{ij}) - \sum_{j=1}^n \sum_{i=1}^{m-1}
-2(x_{ij}) \leq -2n\] and by definition, \[\sum_{j=1}^n \sum_{i=1}^m
\chi(S_i \cap C(K_j)) - \sum_{j=1}^n \sum_{i=1}^{m-1} \chi(F_i \cap C(K_j))
\leq -2n.\] So, \[\chi(S) = \sum_{i=1}^m \chi(S_i) - \sum_{i=1}^{m-1}
\chi(F_i) \leq -2n.\] Whence \[genus(S) \geq n+1\] and \[t(K_1 \#
\dots \# K_n) \geq n.\]
\end{pf}


\begin{thebibliography}{99}

\bibitem{BZ}
{\em Knots}

G. Burde, H. Zieschang,

de Gruyter, Studies in Mathematics 5,
Berlin, New York 

\bibitem{H}
{\em $3$-manifolds}

J. Hempel,

Annals of Math. Studies 86 (1976), Princeton University Press

\bibitem{J}
{\em Lectures on Three-manifold Topology}

W. Jaco,

Regional Conference Series in Mathematics 43 (1981), Amer. Math. Soc.

\bibitem{K}
{\em Problem's in Low-Dimensional Topology}

R. Kirby,
in {\em Geometric Topology -- Proceedings of the 1993 Georgia
International Topology Conference}, ed. by W. H. Kazez, Vol. 2,
Part 2, pp. 35--473

\bibitem{N}
{\em Every two generator knot is prime}

F.H. Norwood,
Proc. Amer. Math. Soc. 86 (1982), 143--147

\bibitem{R}
{\em Knots and Links}

D. Rolfsen,

Publish or Perish, Inc.
Houston, Texas

\bibitem{S1}
{\em Heegaard splittings of compact $3$-manifolds}

M. Scharlemann,

to appear in the {\em Handbook of Geometric Topology}, ed. by R. Daverman
and R. Sherr, Elsevier Press.

\bibitem{ST1}
{\em Thin position for $3$-manifolds}

M. Scharlemann, A. Thompson,
AMS Contemporary Math. 164 (1994) 231--238

\bibitem{Sc}
{\em Additivity of Tunnel Number}

J. Schultens,

preprint

\end{thebibliography}
\end{document}